\documentclass{amsart}

\usepackage{amssymb}
\usepackage{amsrefs}

\newtheorem{thm}{Theorem}
\newtheorem{cor}[thm]{Corollary}
\newtheorem{lemma}{Lemma}
\theoremstyle{remark}
\newtheorem*{remark}{\indent Remark}

\newcommand{\abs}[1]{\left\lvert#1\right\rvert}
\newcommand{\dbars}[1]{\left\lVert#1\right\rVert}

\newcommand{\eulerian}[2]{\genfrac{\langle}{\rangle}{0pt}{}{#1}{#2}}
\DeclareMathOperator{\vol}{Vol} \DeclareMathOperator{\sinc}{sinc}

\begin{document}

\title{Slices, slabs, and sections of the unit hypercube}
\author{Jean-Luc Marichal}
\address{Institute of Mathematics\\
         University of Luxembourg\\
         162A avenue de la Fa\"{\i}encerie\\
         L-1511 Luxembourg\\
         Luxembourg}
\email{jean-luc.marichal@uni.lu}
\author{Michael J. Mossinghoff}
\address{Department of Mathematics\\
         Davidson College\\
         Davidson, North Carolina 28035-6996 USA}
\email{mimossinghoff@davidson.edu}

\date{\today}



\begin{abstract}
Using combinatorial methods, we derive several formulas for the volume of convex bodies obtained by intersecting a unit hypercube with a
half\-space, or with a hyperplane of codimension~1, or with a flat defined by two parallel hyperplanes.
We also describe some of the history of these problems, dating to P\'olya's Ph.D. thesis, and we discuss several applications of these formulas.
\end{abstract}

\maketitle

\section{Introduction}%
\label{sectionIntroduction}

In this note we study the volumes of portions of $n$-dimensional cubes determined by hyperplanes. More precisely, we study \textit{slices}\/
created by intersecting a hypercube with a halfspace, \textit{slabs}\/ formed as the portion of a hypercube lying between two parallel
hyperplanes, and \textit{sections}\/ obtained by intersecting a hypercube with a hyperplane. These objects occur naturally in several fields,
including probability, number theory, geometry, physics, and analysis. In this paper we describe an elementary combinatorial method for
calculating volumes of arbitrary slices, slabs, and sections of a unit cube. We also describe some applications that tie these geometric results
to problems in analysis and combinatorics.

Some of the results we obtain here have in fact appeared earlier in other contexts. However, our approach is entirely combinatorial, in contrast
with most of the other work on this topic. In addition, due to the wide application of cube slices and sections in several fields, it appears
that some of the previous work we describe here (including, for instance, P\'olya's thesis) is not widely cited in the literature, and we hope
that by surveying some of this work here we might make some of the history and prior contributions on these topics better known.

We remark that Zong's recent article \cite{Zong05} and monograph \cite{Zong06} survey many problems and topics regarding $n$-dimensional cubes
and parallelotopes, including questions on sections, projections, triangulations, and inscribed simplices. We refer the reader to these sources
for additional problems and questions about unit cubes.

We employ the following notation throughout. Let $I^n$ denote the $n$-dimensional unit cube placed in the positive orthant of $\mathbb{R}^n$
with one vertex at the origin,
\[
I^n := [0,1]^n.
\]
Also, sometimes it is more convenient to consider the unit cube centered at the origin in $\mathbb{R}^n$, and we let $C^n$ denote this  cube,
\begin{equation}\label{eqnCn}
C^n := \left[-{\textstyle\frac{1}{2}},{\textstyle\frac{1}{2}}\right]^n.
\end{equation}
In particular, it is often more useful to work with $C^n$ rather than $I^n$ in analytic applications, since $C^n$ is simply the ball of radius
$1/2$ with respect to the $\ell_\infty$ norm on $\mathbb{R}^n$.

Next, for a nonzero vector $\mathbf{w}:=(w_1,\ldots,w_n)$ in $\mathbb{R}^n$ and a real number $z$, let $G_{\mathbf{w},z}^n$ denote the halfspace
in $\mathbb{R}^n$ given by
\[
G_{\mathbf{w},z}^{n} := \{\mathbf{x}\in\mathbb{R}^n : \mathbf{w}\cdot\mathbf{x}\leqslant z\},
\]
and let $H_{\mathbf{w},z}^{n-1}$ denote the corresponding hyperplane with normal vector $\mathbf{w}$,
\[
H_{\mathbf{w},z}^{n-1} := \{\mathbf{x}\in\mathbb{R}^n : \mathbf{w}\cdot\mathbf{x}=z\}.
\]
If some component of $\mathbf{w}$ is 0, then clearly any question on slices, slabs, or sections of unit cubes involving $G_{\mathbf{w},z}^{n}$
or $H_{\mathbf{w},z}^{n-1}$ reduces to a problem in a lower dimension, so we assume that each component of $\mathbf{w}$ is nonzero. It is also
clear by symmetry that the volume of a slice, slab, or section is unchanged when the components of $\mathbf{w}$ are permuted.

For a positive integer $n$, let $[n]$ denote the set $\{1,2,\ldots,n\}$, and let $V_n$ denote the set of vertices of $2C^n$, that is,
$V_n:=\{-1,1\}^n$. Also, for $\mathbf{s}\in V_n$, let $\epsilon_{\mathbf{s}}:=\prod_{i=1}^n s_i$. For a real number $p\geqslant 1$, let
$\dbars{\cdot}_p$ denote the usual $\ell_p$ norm on $\mathbb{R}^n$, so $\dbars{\mathbf{w}}_p := (\sum_{i=1}^n \abs{w_i}^p)^{1/p}$, and let
$\dbars{\mathbf{w}}_\infty :=\max\{\abs{w_i}:i\in [n]\}$. We also set $N_{\mathbf{w}}:=\{i\in [n]:w_i<0\}$, and we denote by $A_{\mathbf{w}}$ the
$n\times n$ diagonal matrix whose $i$th diagonal entry is $1$ if $w_i$ is positive and $-1$ if it is negative.
In addition, for any $K\subseteq [n]$, we denote by $\mathbf{1}_K$ the characteristic
vector of $K$ in $\{0,1\}^n$, and for subsets $K_1$ and $K_2$ of $[n]$ we let
$K_1\ominus K_2$ denote their symmetric difference.
Also, for a real number $r$ and nonnegative integer $n$, we set
$r_+^n:= (\max\{r,0\})^n$. Finally, if $\mathbf{x}$ and $\mathbf{y}$ are vectors in $\mathbb{R}^n$, we write $\mathbf{x}\succ\mathbf{y}$ if $x_i>y_i$
for each $i$, and in the same way we write $\mathbf{x}\succcurlyeq \mathbf{y}$ if $x_i\geqslant y_i$ for each $i$.

This paper is organized in the following way. In Section~\ref{sectionSlice} we derive a formula for the volume of an arbitrary slice of a cube
by using an elementary combinatorial argument (see Theorem~\ref{thmSlice}), and we use this to derive some formulas for volumes of central slabs
of cubes (see Theorem~\ref{thmSlab} and Corollary~\ref{corSlab}). We also describe some of the history of these problems here, including
P\'olya's dissertation on volumes of central slabs of cubes. In Section~\ref{sectionSection} we use our formula for slices to determine some
formulas for the volume of a section of a cube (see Theorem~\ref{thmSection} and Corollary~\ref{corSection}), and we describe some of the
history and applications of the problem of bounding this volume. In Section~\ref{sectionApplications} we describe several applications of these
formulas, including some formulas for integrating a polynomial over a cube slice or a cube section, some computations in probability,
and some combinatorial identities, including a geometric interpretation of the Eulerian numbers.

\section{Slices and slabs of cubes}%
\label{sectionSlice}

We derive an exact formula for the volume of a hypercube sliced by a hyperplane of codimension~1. Our method, which generalizes a geometric
approach proposed in 1989 by Denardo and Larraza~\cite{DL}, employs a signed simplicial subdivision of the sliced cube, together with the
inclusion-exclusion principle.

For any $\mathbf{w}\in\mathbb{R}^n$ with $\mathbf{w}\succ\mathbf{0}$, any $z\in\mathbb{R}$, and any $K\subseteq [n]$, define the set
$\Delta_{\mathbf{w},z}^K$ by
\[
\Delta_{\mathbf{w},z}^K := G_{\mathbf{w},z}^n \cap \{\mathbf{x}\in\mathbb{R}^n : \mathbf{x}\succcurlyeq\mathbf{1}_K\}.
\]
This set is empty if $\mathbf{w}\cdot\mathbf{1}_K >z$, and it is an $n$-simplex if $\mathbf{w}\cdot\mathbf{1}_K \leqslant z$. Also, we
immediately see that by translating the set $\Delta_{\mathbf{w},z}^K$ by the vector $-\mathbf{1}_K$ we obtain the set
$\Delta_{\mathbf{w},z-\mathbf{w}\cdot\mathbf{1}_K}^{\varnothing}$.

The following lemma, which is the key result of this paper, yields a signed decomposition of the sliced cube $G_{\mathbf{w},z}^n \cap I^n$ into
the simplices $\Delta_{\mathbf{w},z}^K$.

\begin{lemma}\label{lemmaSieve}
Suppose $\mathbf{w}\in\mathbb{R}^n$ has all positive components, and suppose $z$ is a real number. Then
\begin{equation}\label{eqSieve}
\vol_n\bigl(G_{\mathbf{w},z}^n \cap I^n\bigr) = \sum_{K\subseteq[n]} (-1)^{\abs{K}} \vol_n\bigl(\Delta_{\mathbf{w},z}^K\bigr).
\end{equation}
\end{lemma}

\begin{proof}
By writing $I^n = \left\{\mathbf{x}\in\mathbb{R}^n : \mathbf{x}\succcurlyeq\mathbf{0}\right\} \setminus \bigcup_{i\in[n]}
\left\{\mathbf{x}\in\mathbb{R}^n : \mathbf{x}\succ\mathbf{1}_{\{i\}}\right\}$, we see that
\begin{align*}
\vol_n\bigl(G_{\mathbf{w},z}^n\cap I^n\bigr)
&= \vol_n\biggl(\Delta_{\mathbf{w},z}^{\varnothing} \setminus \bigcup_{i\in [n]} \Delta_{\mathbf{w},z}^{\{i\}}\biggr)\\
&= \vol_n\bigl(\Delta_{\mathbf{w},z}^{\varnothing}\bigr) - \vol_n\biggl(\,\bigcup_{i\in [n]} \Delta_{\mathbf{w},z}^{\{i\}}\biggr),
\end{align*}
since clearly $\Delta_{\mathbf{w},z}^{\{i\}}\subseteq\Delta_{\mathbf{w},z}^{\varnothing}$ for each $i$. Then, by inclusion-exclusion, we find
that
\begin{align*}
\vol_n\biggl(\,\bigcup_{i\in[n]} \Delta_{\mathbf{w},z}^{\{i\}}\biggr)
&= \sum_{\substack{K\subseteq [n]\\K\neq\varnothing}} (-1)^{\abs{K}+1} \vol_n\biggl(\,\bigcap_{i\in K}\Delta_{\mathbf{w},z}^{\{i\}}\biggr)\\
&= \sum_{\substack{K\subseteq [n]\\K\neq\varnothing}} (-1)^{\abs{K}+1} \vol_n\bigl(\Delta_{\mathbf{w},z}^K\bigr),
\end{align*}
which completes the proof.
\end{proof}

We now establish a formula for the volume of an arbitrary slice of a hypercube.

\begin{thm}\label{thmSlice}
Suppose $\mathbf{w}\in\mathbb{R}^n$ has all nonzero components, and suppose $z$ is a real number. Then
\begin{equation}\label{eqnSlice}
\vol_n\bigl(G_{\mathbf{w},z}^n \cap I^n\bigr) = \frac{1}{{n!\prod_{i=1}^n w_i}}\sum_{K\subseteq [n]} (-1)^{\abs{K}} \left(z-\mathbf{w}\cdot\mathbf{1}_K\right)_+^n.
\end{equation}
\end{thm}

\begin{proof}
The result follows from \eqref{eqSieve} if we assume that $\mathbf{w}\succ\mathbf{0}$. In fact, as the volume of the $n$-simplex bounded by the
coordinate hyperplanes and the hyperplane $\mathbf{w}\cdot\mathbf{x}=1$ is $1/(n!\prod_{i=1}^n w_i)$ (see for instance \cite{Ellis}), we
immediately have
\[
\vol_n\bigl(\Delta_{\mathbf{w},z}^K\bigr) = \vol_n\bigl(\Delta_{\mathbf{w},z-\mathbf{w}\cdot\mathbf{1}_K}^{\varnothing}\bigr) =
\frac{\left(z-\mathbf{w}\cdot\mathbf{1}_K\right)^n_+}{n!\prod_{i=1}^n w_i}.
\]
Assume now that $\mathbf{w}\in\mathbb{R}^n$ has all nonzero components. By using the change of variables
$\mathbf{x}'=A_{\mathbf{w}}\mathbf{x}+\mathbf{1}_{N_{\mathbf{w}}}$, it is straightforward to show that
\[
\vol_n\bigl(G_{\mathbf{w},z}^n \cap I^n\bigr) = \vol_n\bigl(G_{A_{\mathbf{w}}\mathbf{w},z-\mathbf{w}\cdot\mathbf{1}_{N_{\mathbf{w}}}}^n\cap
I^n\bigr).
\]
Therefore, since $A_{\mathbf{w}}\mathbf{w}$ has all positive components, we obtain
\[
\vol_n\bigl(G_{\mathbf{w},z}^n \cap I^n\bigr)=\frac{1}{{n!\prod_{i=1}^n w_i}}\sum_{K\subseteq [n]} (-1)^{\abs{K\ominus N_{\mathbf{w}}}}
\left(z-\mathbf{w}\cdot\mathbf{1}_{K\ominus N_{\mathbf{w}}}\right)_+^n,
\]
and the result follows by observing that $\{K\ominus N_{\mathbf{w}}:K\subseteq [n]\}=2^{[n]}$.
\end{proof}

The general formula for the volume of a cube slice seems to have first appeared in a note of Barrow and Smith in 1979 \cite{BS}. The proof there
exploits the simple integration formula for the truncated power function $f(x)=(x-a)^k_+$.
Also, Ueda et al.\ in 1994 \cite{Ueda} and Bradley and Gupta in 2002 \cite{BG} independently derived formulas similar to \eqref{eqnSlice} by analytic means in the context of investigating certain problems in probability (see Section~\ref{subsectionProbability}).
A special case of the formula however occurred
considerably earlier. Certainly from Theorem~\ref{thmSlice} one immediately obtains a formula for the volume of the slab
\[
\{\mathbf{x}\in\mathbb{R}^n : z_1 \leqslant \mathbf{w}\cdot\mathbf{x} \leqslant z_2\} \cap I^n,
\]
for real numbers $z_1$ and $z_2$ with $z_1\leqslant z_2$. In his 1912 dissertation \cite{Polya12}, P\'olya studied the special case of
determining the volume of a central slab of a hypercube, motivated by a question in statistical mechanics. His results are published
in \cite{Polya13}. It is simplest to describe his result here using the centered unit cube $C^n$ from \eqref{eqnCn}. Let
$S_{\mathbf{w},\theta}^n$ denote the central slab with normal vector $\mathbf{w}$ in $\mathbb{R}^n$ and thickness $\theta/\dbars{\mathbf{w}}_2$,
\[
S_{\mathbf{w},\theta}^n := \{\mathbf{x}\in\mathbb{R}^n : \abs{\mathbf{w}\cdot\mathbf{x}}\leqslant\theta/2\},
\]
and assume each component $w_i$ of $\mathbf{w}$ is nonzero.  P\'olya in effect proved that
\begin{equation}\label{eqnPolyaIntegral}
\vol_n\bigl(S_{\mathbf{w},\theta}^n\cap C^n\bigr) = \frac{2}{\pi}\int_0^\infty \frac{\sin\theta x}{x}\,\prod_{i=1}^n \frac{\sin w_i x}{w_i x}\,dx.
\end{equation}
This may be verified by considering the Fourier transform of a convolution of characteristic functions of intervals. He also determined a
formula for this integral by using the residue theorem. We establish his result here as a simple consequence of Theorem~\ref{thmSlice}.

\begin{thm}\label{thmSlab}
Suppose $\mathbf{w}\in\mathbb{R}^n$ has all nonzero components, and suppose $\theta$ is a positive real number. Let $\mathbf{v}$ denote the
vector $(w_1,\ldots,w_n,\theta)$. Then
\begin{equation}\label{eqnPolyaSlice}
\vol_n\bigl(S_{\mathbf{w},\theta}^n\cap C^n\bigr) = \frac{1}{2^n n! \prod_{i=1}^n w_i}\sum_{\mathbf{s}\in V_{n+1}}
\epsilon_{\mathbf{s}}\left(\mathbf{v}\cdot\mathbf{s}\right)^n_+.
\end{equation}
\end{thm}

\begin{proof}
Let $y=\mathbf{w}\cdot\mathbf{1}_{[n]}$. Using Theorem~\ref{thmSlice}, we have
\begin{align*}
\vol_n\bigl(&S^n_{\mathbf{w},\theta}\cap C^n\bigr) = \vol_n\bigl(G^n_{\mathbf{w},y/2+\theta/2}\cap I^n\bigr)
- \vol_n\bigl(G^n_{\mathbf{w},y/2-\theta/2}\cap I^n\bigr)\\
&= \frac{1}{2^n n! \prod_{i=1}^n w_i} \sum_{K\subseteq[n]} (-1)^{\abs{K}}
\bigl( \left(y+\theta-2\mathbf{w}\cdot\mathbf{1}_K\right)_+^n
- \left(y-\theta-2\mathbf{w}\cdot\mathbf{1}_K\right)_+^n \bigr)\\
&= \frac{1}{2^n n! \prod_{i=1}^n w_i} \sum_{\mathbf{s}\in V_{n+1}} \epsilon_{\mathbf{s}}\left(\mathbf{v}\cdot\mathbf{s}\right)_+^n.\qedhere
\end{align*}
\end{proof}

\begin{remark}
In fact, P\'olya considered the dual problem of slicing an $n$-dimensional box of arbitrary dimensions centered at the origin, with a slab
having the fixed normal vector $(1,\ldots,1)$: If $\mathbf{w}\in\mathbb{R}^n$ has all positive components and $P_{\mathbf{w},\theta}^n$ denotes
the polytope determined by the inequalities $\abs{x_i}\leqslant w_i$ for $i\in [n]$ and $\abs{\sum_{i=1}^n x_i}\leqslant\theta$, then P\'olya
showed that
\[
\vol_n\bigl(P_{\mathbf{w},\theta}^n\bigr) = \frac{2^{n+1}}{\pi} \int_0^\infty \frac{\sin\theta x}{x}\prod_{i=1}^n\frac{\sin w_i x}{x}\, dx =
\frac{1}{n!}\sum_{\mathbf{s}\in V_{n+1}} \epsilon_{\mathbf{s}}\left(\mathbf{v}\cdot\mathbf{s}\right)_+^n,
\]
where again $\mathbf{v}=(w_1,\ldots,w_n,\theta)$. We have transcribed his results here for the equivalent problem where the box is fixed and the
normal vector of the slab varies. P\'olya also mentioned that this last formula can be obtained by using a signed simplicial subdivision, as in
the proof of Lemma~\ref{lemmaSieve}.
\end{remark}

P\'olya's formula and its close relatives have been rediscovered several times since then by many people (including the authors, independently).
For example, Borwein and Borwein \cite{BB} proved in 2001 that
\[
\frac{2}{\pi}\int_0^\infty \prod_{i=0}^n \frac{\sin(a_i x)}{x}\,dx = \frac{1}{2^n n!} \sum_{\mathbf{s}\in V_n}
\epsilon_{\mathbf{s}}\alpha_{\mathbf{s}} \bigg(a_0+\sum_{i=1}^n s_i a_i\bigg)^n,
\]
where $a_0$, \ldots, $a_n$ are positive real numbers, and $\alpha_{\mathbf{s}}$ is 1, 0, or $-1$ depending on whether $a_0+\sum_{i=1}^n s_i a_i$
is positive, zero, or negative. It is straightforward to show that this is equivalent to \eqref{eqnPolyaSlice}.

We next derive two additional useful formulations for the volume of a central slab of a cube. Let $V_n^-$ denote the set of $n$-tuples in $V_n$
whose last coordinate is $-1$, and let $V_n^+:=V_n\setminus V_n^-$.

\begin{cor}\label{corSlab}
Suppose $\mathbf{w}\in\mathbb{R}^n$ has all nonzero components, and suppose $\theta$ is a positive real number. If
$\mathbf{v}=(w_1,\ldots,w_n,\theta)$, then
\begin{equation}\label{eqnSlabA}
\vol_n\bigl(S_{\mathbf{w},\theta}^n\cap C^n\bigr) = 1 + \frac{1}{2^{n-1}n!\prod_{i=1}^n w_i}\sum_{\mathbf{s}\in V_{n+1}^-}
\epsilon_{\mathbf{s}} \left(\mathbf{v}\cdot\mathbf{s}\right)^n_+.
\end{equation}
Also, if $\mathbf{v}'=(\theta,w_1,\ldots,w_n)$, then
\begin{equation}\label{eqnSlabB}
\vol_n\bigl(S_{\mathbf{w},\theta}^n\cap C^n\bigr) = \frac{\theta}{w_n} + \frac{1}{2^{n-1}n!\prod_{i=1}^n w_i}\sum_{\mathbf{s}\in V_{n+1}^-}
\epsilon_{\mathbf{s}} \left(\mathbf{v}'\cdot\mathbf{s}\right)^n_+.
\end{equation}
\end{cor}

\begin{proof}
By using the immediate identity
\[
r^n=r_+^n+(-1)^n(-r)^n_+\,,
\]
we easily obtain
\begin{equation}\label{eqn1stDec}
\sum_{\mathbf{s}\in V_{n+1}^+} \epsilon_{\mathbf{s}}\left(\mathbf{v}\cdot\mathbf{s}\right)^n_+ = \sum_{\mathbf{s}\in V_{n+1}^+}
\epsilon_{\mathbf{s}}\left(\mathbf{v}\cdot\mathbf{s}\right)^n + \sum_{\mathbf{s}\in V_{n+1}^-} \epsilon_{\mathbf{s}}\left(\mathbf{v}\cdot\mathbf{s}\right)^n_+.
\end{equation}
Using the multinomial theorem, we also have
\[
\sum_{\mathbf{s}\in V_{n+1}^+} \epsilon_{\mathbf{s}}\left(\mathbf{v}\cdot\mathbf{s}\right)^n=\sum_{\substack{k_1,\ldots,k_{n+1}\geqslant 0\\
k_1+\cdots +k_{n+1}=n}}\frac{n!\,\theta^{k_{n+1}}}{k_1!\cdots k_{n+1}!}\,\prod_{i=1}^n w_i^{k_i}\sum_{\mathbf{s}\in V_{n+1}^+}\prod_{i=1}^n
s_i^{k_i+1}.
\]
The inner sum on the right equals $\prod_{i=1}^n \big((-1)^{k_i+1}+1\big)$, which is~0 for each possible choice of $k_1$, \ldots, $k_{n+1}$
except $k_1=\cdots =k_n=1$ and $k_{n+1}=0$. Therefore,
\begin{equation}\label{eqn2ndDec}
\sum_{\mathbf{s}\in V_{n+1}^+} \epsilon_{\mathbf{s}}\left(\mathbf{v}\cdot\mathbf{s}\right)^n=2^n n!\prod_{i =1}^n w_i,
\end{equation}
and \eqref{eqnSlabA} follows from~\eqref{eqnPolyaSlice}, \eqref{eqn1stDec}, and \eqref{eqn2ndDec}. Equation \eqref{eqnSlabB} may be verified in
the same fashion.
\end{proof}

A formula similar to \eqref{eqnSlabA} also appears in Borwein and Borwein \cite{BB}, where it is used in conjunction with
\eqref{eqnPolyaIntegral} to produce some striking formulas for some integrals involving the $\sinc$ function, where $\sinc(x):=\sin(x)/x$. For
example, since the sum of the reciprocals of the odd primes first exceeds~1 at $p=29$, we see from \eqref{eqnPolyaIntegral} and \eqref{eqnSlabA}
that
\[
\int_0^\infty \sinc(x) \prod_{\substack{p \textrm{\ prime}\\3\leqslant p\leqslant m}} \sinc(x/p)\, dx = \frac{\pi}{2}
\]
for $m\leqslant28$, but
\begin{equation*}
\begin{split}
\int_0^\infty \sinc(x) \prod_{\substack{p \textrm{\ prime}\\3\leqslant p\leqslant29}} \sinc(x/p)\, dx
&= \frac{\pi}{2}\left(1-\frac{54084649^9}{181440\cdot3234846615^8}\right)\\
&=(0.49999999999908993\ldots)\pi.
\end{split}
\end{equation*}
Similar formulas appear in \cite{BBKW} and \cite{BB}.
We mention also that Borwein, Borwein, and Mares \cite{BBM} developed further connections between multivariable $\sinc$ integrals and volumes of
polytopes constructed by intersecting flats formed by pairs of symmetric hyperplanes.
In addition, Barthe and Koldobsky \cite{BK} investigated the problem of minimizing the volume of a central slab of the unit hypercube, if the thickness of the slab is fixed.

The formula \eqref{eqnSlabB} for the volume of a central slab of a cube is
employed in Section~\ref{sectionSection} to obtain a formula for the volume
of a central section of a cube.

\section{Sections of cubes}%
\label{sectionSection}

The convex bodies formed by intersecting a hypercube with a hyperplane are also well-studied in analysis and number theory, especially problems
about bounding their volumes. For example, Vaaler \cite{Vaaler} established a sharp lower bound on the volume of a central slice of $C^n$ by a
$k$-dimensional hyperplane $H_0^k$ that passes through the origin,
\[
\vol_k\bigl(H_0^k \cap C^n\bigr) \geqslant 1,
\]
and Ball \cite{Ball89} determined an upper bound for an arbitrary hyperplane $H^k$ of dimension $k$,
\[
\vol_k\bigl(H^k \cap C^n\bigr) \leqslant 2^{(n-k)/2}.
\]
This is sharp when $k\geqslant n/2$.
Similar problems have been studied in other settings, for example,
bounding the volume of sections of the unit ball of various $\ell_p$ norms in
$\mathbb{R}^n$ \cite{MP},
studying sections of polydiscs in $\mathbb{C}^n$ \cite{OP},
and investigating sections of regular simplices \cite{Webb}.
More information may be found for instance in \cite{Koldobsky}.

The lower bound on the volume of a central section of a hypercube plays an important role in the strengthening of Siegel's lemma by Bombieri and
Vaaler \cite{BV}. The upper bound is of interest for example in the well-known (and now completely resolved) problem of Busemann and Petty:  If
$A$ and $B$ are symmetric convex bodies in $\mathbb{R}^n$, must $\vol_n(A)\geqslant\vol_n(B)$ if $\vol_{n-1}\bigl(H_{\mathbf{u},0}^{n-1}\cap
A\bigr)\geqslant\vol_{n-1}\bigl(H_{\mathbf{u},0}^{n-1}\cap B\bigr)$ for every unit vector $\mathbf{u}$ in $\mathbb{R}^n$? (Recall a set
$A\subseteq\mathbb{R}^n$ is \textit{symmetric}\/ if $-\mathbf{x}\in A$ whenever $\mathbf{x}\in A$.) Ball \cite{Ball88} showed that comparing a
hypercube with a sphere provides a counterexample for dimensions $n\geqslant 10$; more generally, the answer is now known to be affirmative only
for $n\leqslant 4$. See Koldobsky \cite{Koldobsky} or Zong \cite{Zong96} for more information.

Here we obtain some exact formulas for the volume of an arbitrary section of a hypercube. It is possible that more precise information on
volumes of sections of cubes may allow for more exact bounds in applications. Indeed, the precise values of certain central sections of cubes
play an important role in the recent work of Aliev \cite{Aliev}, where in effect the volume of the central section of the cube $C^n$ obtained
with the symmetric hyperplane $H_{\mathbf{\mathbf{1}_{[n]}},0}^{n-1}$ is employed in an improved lower bound on the maximal element of a set of
positive integers having distinct subset sums.

In the following theorem we present one formula for general sections of $I^n$
and another for central sections of $C^n$.
For the sake of clarity we observe that $H^{n-1}_{\mathbf{w},z}$ intersects $C^n$ if and only if $\abs{z} \leqslant \frac 12
\dbars{\mathbf{w}}_1$. Similarly, $H^{n-1}_{\mathbf{w},z}$ intersects $I^n$ if and only if $\abs{z-\frac
12\mathbf{w}\cdot\mathbf{1}_{[n]}}\leqslant \frac 12\dbars{\mathbf{w}}_1$, that is, if and only if
$\mathbf{w}\cdot\mathbf{1}_{N_{\mathbf{w}}} \leqslant z \leqslant \mathbf{w}\cdot\mathbf{1}_{[n] \setminus N_{\mathbf{w}}}$.

\begin{thm}\label{thmSection}
Suppose $\mathbf{w}\in\mathbb{R}^n$ has all nonzero components, and suppose $z$ is a real number. Then
\begin{equation}\label{eqnSectionI}
\vol_{n-1}\bigl(H_{\mathbf{w},z}^{n-1} \cap I^n\bigr) = \frac{\dbars{\mathbf{w}}_2}{{(n-1)!\prod_{i=1}^n w_i}}\sum_{K\subseteq [n]} (-1)^{\abs{K}}
\left(z-\mathbf{w}\cdot\mathbf{1}_K\right)_+^{n-1}.
\end{equation}
In particular, for a central section of $C^n$ we have
\begin{equation}\label{eqnSectionC1}
\vol_{n-1}\bigl(H_{\mathbf{w},0}^{n-1} \cap C^n\bigr) = \frac{\dbars{\mathbf{w}}_2}{{2^{n-1}(n-1)!\prod_{i=1}^n w_i}}\sum_{\mathbf{s}\in V_n}
\epsilon_{\mathbf{s}} \left(\mathbf{w}\cdot\mathbf{s}\right)_+^{n-1}.
\end{equation}
\end{thm}

\begin{proof}
Let $\mathbf{w}\in \mathbb{R}^n$ with all nonzero components, and let $z,\delta\in\mathbb{R}$. Since the distance between the hyperplanes
$H_{\mathbf{w},z+\delta}^{n-1}$ and $H_{\mathbf{w},z}^{n-1}$ is $\abs{\delta}/\dbars{\mathbf{w}}_2$, we have
\begin{align*}
\vol_{n-1}\bigl(H_{\mathbf{w},z}^{n-1} \cap I^n\bigr)
&= \lim_{\delta\to 0} \frac{\vol_n\bigl(G_{\mathbf{w},z+\delta}^n\cap I^n\bigr)-\vol_n\bigl(G_{\mathbf{w},z}^n\cap I^n\bigr)}{\delta/\dbars{\mathbf{w}}_2}\\
&= \dbars{\mathbf{w}}_2\frac{\partial}{\partial z}\vol_n\bigl(G_{\mathbf{w},z}^n\cap I^n\bigr).
\end{align*}
(See for instance Lasserre \cite{Lasserre83}.) Since $\frac{d}{dx}(x^n_+)=n(x^{n-1}_+)$ for any positive integer $n$, then \eqref{eqnSectionI}
follows immediately from \eqref{eqnSlice}. It is then straightforward to check that \eqref{eqnSectionC1} follows by setting
$z=\dbars{\mathbf{w}}_1/2$ in \eqref{eqnSectionI}; it may also be verified by computing $\lim_{\theta\to 0^+}
\dbars{\mathbf{w}}_2\,\frac{f(\theta)}{\theta}=\dbars{\mathbf{w}}_2\, f'(0)$, where $f(\theta):=\vol_n\bigl(S_{\mathbf{w},\theta}^n\cap C^n\bigr)$
from \eqref{eqnPolyaSlice} or \eqref{eqnSlabA}.
\end{proof}

Bradley \cite{Bradley} in effect also computed the volume of an arbitrary section of a hypercube using inclusion-exclusion, obtaining a formula analogous to \eqref{eqnSectionI}.
We think however that the more general treatment developed independently in this paper is of interest as well.

Next, we note an alternative formula for the volume of a central section of a cube, similar to Corollary~\ref{corSlab} for central slabs.

\begin{cor}\label{corSection}
Suppose $\mathbf{w}\in\mathbb{R}^n$ has all nonzero components. Then
\begin{equation}\label{eqnSectionC2}
\vol_{n-1}\bigl(H_{\mathbf{w},0}^{n-1}\cap C^n\bigr) = \frac{\dbars{\mathbf{w}}_2}{w_n} + \frac{\dbars{\mathbf{w}}_2}{2^{n-2}(n-1)!\prod_{i=1}^n w_i}
\sum_{\mathbf{s}\in V_n^-} \epsilon_{\mathbf{s}}\left(\mathbf{w}\cdot\mathbf{s}\right)^{n-1}_+.
\end{equation}
\end{cor}

We omit the proof. The statement \eqref{eqnSectionC2} may be verified by using the method of the proof of Corollary~\ref{corSlab} on the formula
\eqref{eqnSectionC1}, or by computing the derivative of \eqref{eqnSlabB} with respect to $\theta$ at $\theta=0$ and multiplying by
$\dbars{\mathbf{w}}_2$, or by constructing a signed simplicial subdivision of the cube section directly, as in the proof of
Lemma~\ref{lemmaSieve}.

Clearly, when computing the volume of the cube section it is advantageous to permute the components of the normal vector $\mathbf{w}$ first so
that $w_n=\dbars{\mathbf{w}}_\infty$ (negating $\mathbf{w}$ first if necessary so that $w_n>0$) and then use \eqref{eqnSectionC2}, since the number of terms in this sum is only about half the number
required in \eqref{eqnSectionC1}. Also, occasionally the formula for the volume reduces to a very simple expression. For example, we see that
$\vol_{n-1}\bigl(H_{\mathbf{w},0}^{n-1}\cap C^n\bigr)=\dbars{\mathbf{w}}_2/\dbars{\mathbf{w}}_\infty$ precisely when $\dbars{\mathbf{w}}_\infty
\geqslant \frac{1}{2}\dbars{\mathbf{w}}_1$; in this case the cube section is simply a parallelotope. Also, it is straightforward to verify that
the expression for the volume has just two terms when $\dbars{\mathbf{w}}_\infty < \frac{1}{2}\dbars{\mathbf{w}}_1 \leqslant
\dbars{\mathbf{w}}_\infty + \min\{\abs{w_i} : i\in [n]\}$; here the cube section is a parallelotope with two opposite simplices shorn off. Of
course, similar simplifications occur in analogous situations for the cube slab formulas \eqref{eqnSlabA} and \eqref{eqnSlabB}.

\section{Applications}%
\label{sectionApplications}

We briefly describe four applications for these formulas for the volumes of slices, slabs, and sections of hypercubes: obtaining some formulas
for integrating functions over cube slices and sections, performing some calculations in probability, investigating a geometric interpretation
of the Eulerian numbers, and deriving a particular combinatorial identity.

\subsection{Integrating over cube slices and
sections}\label{subsectionIntegrate}

For a function $f:\mathbb{R}^n\to\mathbb{R}$, the change of variables
$\mathbf{x}'=A_{\mathbf{w}}\mathbf{x}+\mathbf{1}_{N_{\mathbf{w}}}-\mathbf{1}_{K}$, combined with the signed simplicial decomposition
\eqref{eqSieve}, immediately provides the integration formula
\[
\int_{G_{\mathbf{w},z}^n \cap I^n}f(\mathbf{x})\, d\mathbf{x} = \sum_{K\subseteq [n]}
(-1)^{\abs{K}}\int_{\Delta_{A_{\mathbf{w}}\mathbf{w},z-\mathbf{w}\cdot\mathbf{1}_{K\ominus
N_{\mathbf{w}}}}^{\varnothing}}f(A_{\mathbf{w}}\mathbf{x}'+\mathbf{1}_{K\ominus N_{\mathbf{w}}})\, d\mathbf{x}',
\]
or, equivalently,
\begin{equation}\label{eqIntForm}
\int_{G_{\mathbf{w},z}^n \cap I^n}f(\mathbf{x})\, d\mathbf{x} = \sum_{K\subseteq [n]}
(-1)^{\abs{K}+\abs{N_{\mathbf{w}}}}\int_{\Delta_{A_{\mathbf{w}}\mathbf{w},z-\mathbf{w}\cdot\mathbf{1}_{K}}^{\varnothing}}
f(A_{\mathbf{w}}\mathbf{x}'+\mathbf{1}_{K})\, d\mathbf{x}',
\end{equation}
so integrating a function over the sliced cube reduces to integrating it over standard $n$-simplices. Moreover, if $f$ is continuously
differentiable, we also have (see Lasserre \cite{Lasserre98}*{Lemma 2.2})
\[
\int_{H_{\mathbf{w},z}^{n-1} \cap I^n}f(\mathbf{x})\, d\mathbf{x} = \dbars{\mathbf{w}}_2\frac{\partial}{\partial z}\,\int_{G_{\mathbf{w},z}^n
\cap I^n}f(\mathbf{x})\, d\mathbf{x}.
\]
The case of polynomials in $n$ variables is particularly interesting, since these functions can be integrated over simplices with exact formulas
(see \cite{Lasserre01} or \cite{Lasserre98}).
For instance, if $\mathbf{w}\succ\mathbf{0}$ and $\alpha_1$, \ldots, $\alpha_n$ are positive integers then one may easily show that
\[
\int_{\Delta_{\mathbf{w},z-\mathbf{w}\cdot\mathbf{1}_K}^{\varnothing}}\prod_{i=1}^n x_i^{\alpha_i}\, d\mathbf{x} =
\frac{\left(z-\mathbf{w}\cdot\mathbf{1}_K\right)_+^{n+\sum_i\alpha_i}\, \prod_{i=1}^n\alpha_i!}{\bigl(n+\sum_i\alpha_i\bigr)!\,
\prod_{i=1}^nw_i^{\alpha_i+1}}.
\]
Such formulas are useful, for instance, in some computations in probability, as we describe next.

\subsection{Calculations in probability}\label{subsectionProbability}

As Barrow and Smith \cite{BS} indicated, the volume of the sliced unit cube has the following probabilistic interpretation: Suppose $X_1$,
\ldots, $X_n$ are random variables which are independent and uniformly distributed in $[0,1]$. Then the volume formula \eqref{eqnSlice} yields
the cumulative distribution function of the random variable $Y_{\mathbf{w}}=\sum_{i=1}^n w_iX_i$. It follows that, up to the factor
$\dbars{\mathbf{w}}_2$, formula \eqref{eqnSectionI} yields the corresponding probability density function.

We also remark that, by using three different methods (Fourier transform, operational calculus, and integration by parts), Ueda et
al.\ \cite{Ueda} derived an explicit formula for the distribution of the sum of independent random variables $X_1$, \ldots, $X_n$, where $X_i$ is uniformly distributed on the real interval $[-a_i,a_i]$.
The distribution function is exactly $\vol_n\bigl(G_{\mathbf{w},z}^n \cap C^n\bigr)$, with each $w_i=2a_i$. Later Bradley and Gupta \cite{BG}, apparently
unaware of this formula, rediscovered the corresponding density \eqref{eqnSectionI} by means of the Fourier transform, combined with geometric
considerations. Prior contributions on special cases are also discussed in \cite{BG}, \cite{Chew}, and \cite{Ueda}, and the references therein.

Finally, it is noteworthy that the integration formula \eqref{eqIntForm} makes it possible to consider other probability distributions. Let
$f:\mathbb{R}^n\to\mathbb{R}$ be the probability density function of a set of continuous random variables $X_1$, \ldots, $X_n$, and assume that
$f$ is supported on $I^n$, that is, $f(\mathbf{x})=0$ if $\mathbf{x}\notin I^n$. Then the cumulative distribution function of the linear
combination $Y_{\mathbf{w}}=\sum_{i=1}^n w_i X_i$ of these random variables is
\begin{equation}\label{eqnCumulDist}
\Pr[Y_{\mathbf{w}} \leqslant z] = \int_{G_{\mathbf{w},z}^n \cap I^n}f(\mathbf{x})\, d\mathbf{x}.
\end{equation}
For instance, for independent beta variables, the density function is given by
\[
f(\mathbf{x})=\prod_{i=1}^n \frac{x_i^{\alpha_i-1}(1-x_i)^{\beta_i-1}}{B(\alpha_i,\beta_i)},
\]
where each $\alpha_i$ and $\beta_i$ is a positive real number and $B(x,y)$ denotes the usual beta function,
$B(x,y):=\Gamma(x)\Gamma(y)/\Gamma(x+y)$. If each $\alpha_i$ and $\beta_i$ is in fact a positive integer, then $f$ is a  polynomial in $n$
variables, and the integral \eqref{eqnCumulDist} can therefore be computed exactly with \eqref{eqIntForm}.

\subsection{Eulerian numbers}\label{subsectionEulerian}

Recall that the \textit{Eulerian number}\/ $\eulerian{n}{k}$ records the number of orderings $a_1$, \ldots, $a_n$ of $[n]$ having exactly $k$
ascents, so $a_i<a_{i+1}$ for precisely $k$ values of $i$. Let $\Xi_k^n$ denote the portion of the hypercube $I^n$ lying between the hyperplanes
$\sum_{i=1}^n x_i = k$ and $\sum_{i=1}^n x_i = k+1$,
\[
\Xi_k^n :=  \{\mathbf{x}\in\mathbb{R}^n : k \leqslant \mathbf{x} \cdot \mathbf{1}_{[n]} \leqslant k+1\}\cap I^n.
\]
Laplace seems to have first investigated this problem \cite{Laplace}*{pp.\ 257--260} in the context of a question in probability, implicitly showing that the volume of $\Xi_k^n$ is $\frac{1}{n!}\eulerian{n}{k}$.
Stanley \cite{Stanley} supplied an analytic proof of this fact, answering a question of Foata, by exhibiting a measure-preserving transformation of $I^n$ that maps
$\Xi_k^n$ to the set $\{\mathbf{x}\in I^n : x_i<x_{i+1} \textrm{\ for exactly $k$ values of $i$}\}$, which clearly has volume $\frac{1}{n!}\eulerian{n}{k}$.
Schmidt and Simion \cite{SS} found a similar proof, employing a variation of Stanley's map to show that $\Xi_k^n$ may be partitioned into $\eulerian{n}{k}$ simplices, each with volume $1/n!$.
Chakerian and Logothetti \cite{CL} also investigated
combinatorial properties exhibited by certain slabs of unit cubes, and
discussed the appearance of the Eulerian numbers in this context.

Our proof of Theorem~\ref{thmSlice}, together with the following calculation, provides an alternative combinatorial derivation for the volume of $\Xi_k^n$.
We compute
\begin{align*}
\vol_n\bigl(\Xi_k^n\bigr) &= \vol_n\bigl(G_{\mathbf{1}_{[n]},k+1}^n\cap I^n\bigr) -
\vol_n\bigl(G_{\mathbf{1}_{[n]},k}^n\cap I^n\bigr)\\
&= \frac{1}{n!}\biggl(\,\sum_{j=0}^{k+1} (-1)^j\binom{n}{j}(k+1-j)^n -
\sum_{j=0}^k (-1)^j\binom{n}{j}(k-j)^n \biggr)\\
&= \frac{1}{n!} \sum_{j=0}^{k+1} (-1)^j \binom{n+1}{j}(k+1-j)^n\\
&= \frac{1}{n!}\eulerian{n}{k},
\end{align*}
using a standard identity for Eulerian numbers (see for instance \cite{Comtet}*{section 6.5}).

\subsection{A combinatorial identity}\label{subsectionCombIdent}

Let $\lambda$ be a nonnegative real number, and let $\mathbf{w}\in\mathbb{R}^n$ have all positive components. Setting
$z=\lambda+\mathbf{w}\cdot\mathbf{1}_{[n]}$, then certainly $\vol_n\left(G_{\mathbf{w},z}^n\cap I^n\right)=1$, and using \eqref{eqnSlice} we
find a simple geometric interpretation of the following combinatorial identity for the case $p=n$:
\begin{equation}\label{eqnCombId}
\sum_{K\subseteq [n]}(-1)^{\abs{K}}\,\Bigl(\lambda + \sum_{i\in K} w_i\Bigr)^p =
\begin{cases}
0, & \mbox{if $0\leqslant p< n$,}\\
(-1)^n n!\ \prod_{i=1}^n w_i, & \mbox{if $p=n$.}
\end{cases}
\end{equation}
The case $0\leqslant p<n$ follows by differentiating \eqref{eqnSlice} with respect to $\lambda$, since the volume is constant for
$\lambda\geqslant0$.
A similar observation is made by Bradley \cite{Bradley}.
Further, since \eqref{eqnCombId} is a polynomial identity, it follows that it is in
fact valid for any $\mathbf{w}\in\mathbb{C}^n$ and any $\lambda\in\mathbb{C}$.
This identity may also be verified by using the multinomial theorem, or by using a generating function as in \cite{BB} and (for the case $\lambda=0$) \cite{Ueda}.

\end{document}